\documentclass[12pt] {article}
\usepackage[cp1251]{inputenc}
\usepackage[english, russian] {babel}
\usepackage {amsmath}
\usepackage {amsthm, amssymb}

\usepackage[dvips]{graphicx}
\usepackage {graphicx}
\newcommand {\ib}[1] {\textit{\textbf{#1}}}
\newcommand{\aff}[1] {\mathop{\rm aff}\nolimits\left({#1}\right)}
\newcommand{\lin}[1]{\mathop{lin}\left(#1\right)}

\newcommand{\epi}{\mathop{epi}}
\newcommand{\cone}{\mathop{cone}}
\newcommand{\conv}[1]{\mathop{conv}\left(#1\right)}
\newcommand {\polar}[1] {{#1}^{\triangle}}
\newcommand {\dubpolar}[1] {{#1}^{\triangle\triangle}}
\newcommand {\polface}[1] {{#1}^{\diamond}}
\newcommand {\dubpolface}[1] {{#1}^{\diamond\diamond}}
\newcommand {\normalfan}[1] {\mathcal N\left(#1\right)}

\renewcommand{\int}[1]{\mathop{Int}\left(#1\right)}

\newcommand{\n}{\mathbf{n}}

\newcommand{\braces}[1]{\left\{#1\right\}}

\newtheorem{theorem}{Теорема}[section]

\theoremstyle{plain}

\newtheorem{proposition}[theorem]{Предложение}
\newtheorem{lemma}[theorem]{Лемма}
\newtheorem{corollary}[theorem]{Следствие}
\newtheorem{definition}[theorem]{Определение}

\begin{document}

\title{Тесные веера и их канонические нормировки}
\author{А.А.~Гаврилюк}


\date{\today}

\begin{abstract}
  В статье вводится новый класс полных вееров в евклидовом пространстве --- тесные веера. Такие
	веера определяются через свойство локальной симметрии разбиения в некоторой грани, которое 
	определено для произвольных разбиений евклидова пространства на полиэдры.
	
	Тесные веера оказываются связаны с теорией параллелоэдров. Доказана теорема о том, что веер 
	схождений параллелоэдров разбиения в некоторой грани является тесным. При помощи этого свойства
	предложено новое доказательство теоремы Делоне о 5 типах схождений параллелоэдров в грани 
	коразмерности 3. 
	
	Канонические нормировки разбиений евклидова пространства --- это специальные положительно 
	определённые функции, заданные на гипергранях разбиения. Известно, что существование таких 
	функций связано с возможностью представить разбиение в качестве проекции некоторой полиэдральной 
	поверхности. В работе доказано, что полный веер в евклидовом пространстве имеет каноническую 
	нормировку тогда и только тогда, когда его можно представить как веер граней некоторого 
	многогранника. 
	
	В качестве приложения предложенных техник, доказано, что звёзды граней коразмерности 2 и 3 в
	разбиении на параллелоэдры имеют каноническую нормировку.
		
  \bigskip
  \noindent {\bf Ключевые слова:} тесный веер, политопальность, каноническая нормировка, 
	  параллелоэдры, женератриса
  
  \bigskip
  \noindent УДК 514.174 + 514.87
\end{abstract}

\maketitle

\tableofcontents
\newpage

\section{Введение, основные понятия и конструкции}\label{sectionIntro}
\subsection{Многогранные веера. Веера граней и веера схождений}

  Важный вопрос теории разбиений евклидовых пространств --- существует ли женератриса данного 
	разбиения. Этот вопрос оказывается тесно связан с понятием канонической нормировки разбиения. 
	В работе \cite{Gavriliuk_liftings} доказано, что женератриса нормального разбиения существует
	тогда и только тогда, когда существует каноническая нормировка данного разбиения.

	\begin{definition}\label{definitionIntroPolyhedron}
	  Под \emph{полиэдром} в данной работе понимается множество точек в евклидовом пространстве
		$\mathbb E^d$, которое представляется в виде пересечения не более чем счётного набора замкнутых
		полупространств $H^+_i$, ограниченных гиперплоскостями $H_i$, и для которого любой шар в
		$\mathbb E^d$ пересекает лишь конечное подмножество из $\{H_i\}$.
		
		Гиперплоскость $H$ называется \emph{опорной} для полиэдра $\mathcal P$, если $H\bigcap \mathcal 
		  P = H\bigcap \partial\mathcal P \neq \emptyset$.
		
		\emph{Гранью} полиэдра называется его пересечение с произвольной его опорной гиперплоскостью.
	\end{definition}
	
	В случае ограниченного множества $\mathcal P$ множество $\{H_i\}$ конечно, и данное определение 
	совпадает со стандартным определением выпуклого многогранника, его граней и опорных 
	гиперплоскостей \cite{Ziegler}.
	
	Если полиэдр $\mathcal P$ неограничен, то легко видеть, что он является замкнутым выпуклым 
	множеством (в вместе с любой парой точек $u, v \in \mathcal P$ содержит также весь отрезок $uv$).
	Ограниченные грани полиэдра являются ограниченными выпуклыми многогранниками, а неограниченные 
	--- полиэдрами меньшей размерности. 
	
	Отсюда легко выводится, что $\mathcal P$ имеет не более, чем счётное число граней (произвольной 
	размерности). Также $\mathcal P$ локально конечен --- любой шар в $\mathbb E^d$ пересекает 
	конечное множество граней из $\mathcal P$.
	
	\begin{definition}
  	\emph{Разбиением} $d$-мерного евклидова пространства на полиэдры называется набор замкнутых
	  $d$-мерных полиэдров $\mathcal T = \{C_1, C_2, \ldots\}$ в $\mathbb E^d$, элементы которого 
		не пересекаются по внутренним точкам и в совокупности покрывают всё пространство: 
	  $\bigcup\limits^{\infty}_{i=1}C_i = \mathbb E^d$. 
		
		Полиэдры, составляющие такой набор, называются \emph{ячейками} разбиения. Грань $F$ 
		размерности $k \leqslant d$ считается принадлежащей разбиению $\mathcal T$, если она является 
		гранью некоторой ячейки разбиения.
	\end{definition}
	
	\begin{definition}
  	\emph{Нормальным разбиением} называется разбиение ``грань-в-грань'', то есть в котором 
		пересечение любых двух ячеек либо пусто, либо является гранью каждого из этих полиэдров 
		(возможно, самим полиэдром). 
	\end{definition}
	
	\begin{definition}
	  \emph{Полиэдральный комплекс} $\mathcal C$ --- это такой набор полиэдров в $\mathbb E^d$,
		что\\
		\noindent 1)\quad пустой полиэдр принадлежит $\mathcal C$;\\
		\noindent 2)\quad если $P \in \mathcal C$, то все грани полиэдра $P$ тоже принадлежат 
		  $\mathcal C$;\\
		\noindent 3)\quad пересечение $P\bigcap Q$ двух полиэдров $P, Q \in \mathcal C$ является 
		  гранью как полиэдра $P$, так и полиэдра $Q$.
			
		\emph{Размерностью} полиэдрального комплекса $\mathcal C$ называется наибольшая размерность
		  полиэдра из $\mathcal C$.
			
		Элементы комплекса $\mathcal C$ называются его \emph{гранями}.
	\end{definition}
	
	Множество всех граней ячеек, входящих в нормальное разбиение, образует полиэдральный комплекс 
	\cite{Ziegler}. 
	
		\begin{definition}
  	Грани размерности $(n-1)$ в комплексе размерности $n$ называются гипергранями, а грани 
		размерности $(n-2)$ такого комплекса называются \emph{корёбрами}.
	\end{definition}
	
	\begin{definition}
	  \emph{Звездой} $k$-мерной грани $F$ нормального разбиения $\mathcal T$ называется множество всех
		граней в $\mathcal T$, содержащих $F$, включая саму грань $F$ и ячейки разбиения.
		
		Про ячейки $C_1, C_2, \ldots, C_k$ нормального разбиения $\mathcal T$ говорят, что они 
		\emph{сходятся} в $k$-грани $F \subset \mathcal T$, если они представляют все ячейки из звезды 
		грани $F$.
		
		Звезда грани $F$ обозначается $St(F)$.
	\end{definition}
	
	\begin{definition}
  	\emph{Полиэдральным конусом} называется множество всех неотрицательных линейных комбинаций 
	  $\{\mathbf{p} + \alpha_1\mathbf{v}_1 + \ldots + \alpha_k\mathbf{v}_k\}, \alpha_i \geqslant 0$ для 
	  конечного набора векторов $\{\mathbf{v}_i\}$ и точки $\mathbf{p}$ в $\mathbb R^d$ 
	  \cite{Ziegler}.
	\end{definition}
	
	\begin{definition}
  	\emph{Веером} $\mathcal F$ в $\mathbb E^d$ называется семейство полиэдральных конусов 
	  $\{C_1, C_2, \ldots C_n\}$, которое является полиэдральным комплексом и пересечение всего
  	семейства непусто. То есть\\
	  $\bullet$ всякая грань конуса $C_i \in \mathcal F$ также является конусом из $\mathcal F$,\\
  	$\bullet$ $C_i \bigcap C_j$ является непустой гранью каждого из $C_i$ и $C_j$.
	\end{definition}

  \begin{definition}	
  	Веер $\mathcal F$ называется \emph{полным}, если конусы из $\mathcal F$ покрывают всё 
	  пространство $\mathbb E^d$, то есть если $\mathcal F$ является разбиением. 
	\end{definition}
	
	Все веера, которые рассматриваются далее, являются полными.

  \begin{definition}	
  	Легко видеть, что в каждом веере есть минимальная по включению непустая грань, которая 
	  принадлежит всем остальным конусам веера. Если эта грань --- вершина, то конус называется 
  	\emph{заострённым}.
	\end{definition}
	
	Веера моделируют локальную геометрию разбиений на многогранники в окрестностях отдельных граней. 
	В работах \cite{Delone} и \cite{Magazinov3Tight} доказано, что для граней коразмерности 3 в 
	разбиении на $d$-мерные параллелоэдры существует 5 различных комбинаторных типов схождений ячеек,
	то есть им соответствует 5 различных типов трёхмерных вееров. Новое доказательство этого факта
	предложено в данной работе.
	
		\begin{definition}\label{definitionIntroLifting}
  	\emph{Подъёмом} $k$-мерного полиэдра $C$ из разбиения $\mathcal T \subset \mathbb E^d$ 
		называется такой $k$-мерный полиэдр $\widetilde{C}$ в расширенном пространстве $\mathbb 
		  E^{d+1}$, содержащем $\mathbb E^d \supset \mathcal T$ как гиперплоскость, который 
		переходит в $C$ при ортогональной проекции на эту гиперплоскость. Также подъёмом называется 
		процедура построения такого полиэдра $\widetilde{C}$ по полиэдру $C$.
	\end{definition}
	
	Легко проверить, что $\aff{\widetilde{C}}$ не перпендикулярно данному $\mathbb E^d$. Отсюда 
	вытекает, что при ортогональной проекции также сохраняется размерность произвольной 
	$l$-грани $\widetilde{F} \subset \widetilde{C}$, и что существует соответствующая $l$-грань 
	$F \subset C$, для которой $\widetilde{F}$ является подъёмом. 
	
	\begin{definition}
	  Если полиэдр $P$ является гранью полиэдра $Q$, то это обозначается $P \prec Q$.
		
		Две грани $P$ и $Q$ полиэдрального комплекса в $\mathbb E^d$ называются \emph{инцидентными},
		если одна из них является гранью другого (возможно, несобственной), то есть если $P \prec Q$ 
		или $Q \prec P$.
	\end{definition}
	
	Таким образом, ортогональная проекция $\widetilde{C}$ на $C \subset\mathbb E^d$ сохраняет 
	инцидентности граней полиэдра $\widetilde{C}$, сохраняя при этом порядок включения граней. 
		
	\begin{definition}
  	\emph{Строго выпуклым} в теории многогранников называется такой полиэдр, что любые три точки 
		его границы, лежащие на одной прямой, принадлежат некоторой собственной грани этого полиэдра.
	\end{definition}

  \begin{definition}\label{definitionIntroGeneratix}
  	\emph{Женератрисой} $\mathcal G_{\mathcal T}$ нормального локально-конечного разбиения $\mathcal 
	    T \subset \mathbb E^d$ на выпуклые многогранники называется полиэдральная $d$-мерная
  	поверхность в расширенном пространстве $\mathbb E^{d+1} \supset \mathbb E^d \supset \mathcal T$ 
	  со следующими свойствами:\\
  	1. \quad $\mathcal G_{\mathcal T}$ является границей строго выпуклого $(d+1)$-мерного полиэдра в 
	    $\mathbb E^{d+1}$.\\
  	2. \quad Ортогональная проекция $\mathcal G_{\mathcal T}$ на $\mathbb E^d \supset \mathcal T$ ---
	    однозначное отображение, образ которого покрывает всё $\mathbb E^d$.\\
  	3. \quad Ортогональная проекция переводит произвольную грань из $\mathcal G_{\mathcal T}$ в грань
	    той же размерности из $\mathcal T$, то есть задаёт изоморфим $\mathcal G_{\mathcal T}$ и 
		  $\mathcal T$ как полиэдральных комплексов.
		
	  Женератрису разбиения мы также будем называть \emph{подъёмом разбиения}, так как каждая грань
	  женератрисы является подъёмом соответствующей грани из разбиения.
	
  	Очевидно, что $\mathcal G_{\mathcal T}$ является графиком некоторой	выпуклой кусочно-линейной 
	  функции $G(x): \mathbb E^d \to \mathbb R$, которая линейна на каждой ячейке из $\mathcal T$. 
  	Строго выпуклым полиэдром из определения женератрисы является надграфик $\braces{(x, y) \in 
		  \mathbb E^d \times \mathbb R|~y \geqslant G(x)} \subset \mathbb E^{d+1}$ функции $G(x)$. 
		Такую функцию $G(x)$ мы также будем называть женератрисой разбиения $\mathcal T$. Надграфик 		 
		функции $G(x)$ обозначается $\epi G$. 
	\end{definition}
	
	Г.Ф. Вороной использовал подъёмы специальных разбиений на параллелоэдры до описанной около 
	некоторого параболоида поверхности. Эта идея была переоткрыта в \cite{Edelsbrunner} и 
	используется в приложениях как алгоритм построения диаграмм Вороного и разбиений Делоне для 
	произвольных локально-конечных точечных множеств. 
	
	В работе \cite{Davis} показано, что при $d \geqslant 3$ для произвольного конечного нормального 
	разбиения на примитивные многогранники существует женератриса. Показано, что для таких разбиений 
	женератрису можно построить, поднимая ячейки последовательно одна за одной.
	
	В работах \cite{Aurenhammer, McMullen} показано, что женератриса разбиения на выпуклые 
	многогранники существует тогда и только тогда, когда для этого разбиения существует 
	\emph{ортогонально-дуальное} множество.
	
	Пусть задано разбиение $\mathcal T$ пространства $\mathbb E^d$ и $F$ --- коребро в нём. При 
	проекции вдоль $\aff{F}$ на двумерную плоскость $\alpha$ дополнительную к $\aff{F}$ грань $F$
	перейдёт в точку $f$. Ячейки из $\mathcal T$, содержащие $F$, перейдут в двумерные многоугольники,
	сходящиеся в $f$, гиперграни, содержащие $F$, перейдут в одномерные рёбра между некоторыми парами 
	полученных многоугольников.
	
	\begin{definition}\label{definitionIntroCycling}
	  \emph{Циклическим обходом ячеек} из разбиения $\mathcal T$ \emph{вокруг коребра $F$} назовём 
		такую нумерацию $C_1, C_2, \ldots, C_k$ ячеек, сходящихся в $F$, что их проекции в $\alpha$ 
		следуют в одном направлении обхода вокруг точки $f$.
		
		Аналогично определяется \emph{циклический обход гиперграней} $F_1, \ldots, F_k$ из разбиения 
		$\mathcal T$ вокруг коребра $F$.
		
		Очевидно, циклические обходы вокруг коребра делятся на две группы циклически эквивалентных 
		между собой. Эти две группы называются \emph{направлениями циклического обхода} вокруг $F$. 
		
		С точностью до циклического сдвига, произвольные два обхода из данных двух групп противоположны 
		друг другу. Два различных направления циклического обхода называются \emph{противоположными}.
	\end{definition}
	
	\begin{definition}\label{definitionIntroNeighbourly}
	  Две $k$-мерные грани $n$-мерного комплекса называются \emph{смежными} если пересекаются
		друг с другом по грани размерности $(k-1)$.
	\end{definition}

	В \cite{RyshRyb} доказано, что для разбиений, \emph{заданных канонически}, существует женератриса.
	Канонически заданное разбиение --- это такой способ записи уравнений $E_{i, i+1}(x) = 0$ для 
	гиперплоскостей, разделяющих смежные ячейки $P_i$ и $P_{i+1}$, что сумма $\sum E_{i, i+1}(x)$ 
	уравнений при циклическом обходе вокруг произвольной $(d-2)$-грани $F$ тождественно равна нулю. 

  Единичную нормаль к гиперграни можно выбрать двумя противоположными способами.
  \begin{definition}
	  Выбор единичных нормалей $\mathbf{n}_1, \ldots, \mathbf{n}_k$ к гиперграням $F_1, \ldots, F_k$, 
		содержащим коребро $F$, называется \emph{согласованным с циклическим обходом} ячеек $C_1, 
		  \ldots, C_k$, сходящихся в $F$, если каждая нормаль $\mathbf{n}_i$ к гиперграни $C_i \bigcap 
			C_{i+1}$ указывает направление от $C_i$ к $C_{i+1}$ (считая $C_{k+1} ~:= C_1$).
	\end{definition}

	\begin{definition}\label{definitionIntroTorsion}
  	Пусть на гипергранях разбиения $\mathcal T$ задана положительная функция (нормировка) 
	  $s:~{\mathcal F}^{d-1}~\to~\mathbb R_+$. Пусть $F_1, F_2, \ldots, F_k$ --- циклический обход 
		гиперграней разбиения, которые содержат коребро $F$. Пусть $\mathbf{n}_i, \ldots, \mathbf{n}_k$ 
		--- согласованный выбор единичных нормалей к этим гиперграням. {\it Кручением} $\Delta_s$ 
		нормировки $s$ вокруг коребра $F$ называется величина 
  	$$\Delta_s(F) := \sum^k_{i=0} s(F_i)\mathbf{n}_i$$
	\end{definition}
	
	Кручение определено с точностью до знака и, вообще говоря, зависит от	выбора направления обхода. 
	Однако нас будет интересовать лишь равенство или неравенство кручения нулю, что от направления 
	обхода не зависит. 

  \begin{definition}
  	Нормировка $s$ множества гиперграней $\mathcal F^{d-1}$ локально-конечного разбиения 
		$\mathcal T$ на выпуклые многогранники называется \emph{канонической}, если кручение $s$ вокруг 
		произвольного коребра $F \subset \mathcal T$ равно нулю.
	\end{definition}
	
	Канонические нормировки разбиения обобщают каноническое задание разбиения из работы 
	\cite{RyshRyb}. Левая часть равенства, определяющего каноническое задание, имеет вид $x\cdot 
	  \Delta(F) + \sum d_{i, i+1}$, где $\Delta$ является кручением некоторой нормировки, а $\sum 
		d_{i, i+1}$ --- это сумма свободных коэффициентов соответствующих уравнений. Таким образом, для 
	того, чтобы задать разбиение канонически, необходимы каноническая нормировка и согласованный 
	набор свободных коэффициентов. 
	
	Задача построения канонической нормировки разбиения сводится к двум другим: задаче построения 
	множества локальных канонических нормировок и задаче ``склейки'' этих нормировок в одну глобальную. 
	
	В данной работе строится новый метод построения локальных нормировок при помощи политопальных 
	вееров. Рассматривается новый объект --- \emph{тесный веер}, который обобщает понятие схождения
	в грани разбиения на параллелоэдры. Доказываются теоремы, что	в малых размерностях (2 и 3) эти 
	веера политопальны и порождают локальные нормировки. Доказанные результаты теории тесных вееров
	применяются для нового доказательства классификации Делоне схождения параллелоэдров в гранях
	коразмерности 3 и существования локальной канонической нормировки звёзд этих граней.
	
	В данной главе используются три стандартные конструкции, определяющие полные веера --- это 
	конструкции \emph{веера граней, веера схождений в грани разбиения} и конструкция 
	\emph{нормального веера} заданного многогранника. 
	
	\begin{definition}
  	\emph{Конической оболочкой} (с вершиной $v \in \mathbb E^d$) множества $M \subset \mathbb E^d$ 
		называется множество неотрицательных конечных линейных комбинаций 
		$\cone_v(M) = \left\{v + \sum\limits^k_{i = 1}\alpha_i\overrightarrow{vp_i}
	  |~ \alpha_i \geqslant 0, p_i \in M\right\}$. 
	\end{definition}
	
	Из определения следует, что коническая оболочка $M$ --- выпуклое множество. Если $M$ многогранник, 
	то $\cone_v(M)$ --- многогранный конус.
	
	\begin{definition}
	  \emph{Веером граней} $\mathcal F_{v, M}$ выпуклого $d$-мерного многогранника $M \subset \mathbb 
		  E^d$ называется набор конических оболочек с вершиной $v \in \int{M}$ граней многогранника $M$:
		$$\mathcal F_{v, M} = \left\{\cone_{v}(F)|~ F \hbox{ --- грань } M\right\}$$
		
		Веер, который можно представить как веер граней, называется \emph{политопальным} (polytopical).
	\end{definition}
	
	Очевидно, что веер граней является полным заострённым веером. 
	
	Возникает	естественное соответствие: произвольная $k$-грань $F \subset M$ ($0 \leqslant k 
	  \leqslant d-1$) соответствует $(k+1)$-мерному конусу $\widetilde{F} = \cone_v(F)$. По 
	определению, пустая грань $\emptyset$ многогранника $M$ соответствует вершине $v$ данного веера, 
	то есть $\widetilde{\emptyset} = v$. 
	
	Очевидно, что данное соответствие сохраняет инцидентности граней. В частности, непересекающиеся 
	грани соответствуют конусам, пересекающимся по вершине $v$.
	
	Легко видеть, что для всякого выпуклого многогранника $M$ существует много вееров граней. Но не 
	всякий веер можно представить как веер граней. Примеры непредставимых вееров можно найти в
	книгах \cite{Ziegler, Fulton, BuhshtaberPanov}. 
		
	Эти примеры существуют уже в размерности 3. В размерности 2 для произвольного полного 
	заострённого веера легко построить порождающий многоугольник, откуда следует известное 
	предложение:
	
	\begin{proposition}\label{propositionTightFans2dPolytopal}~{
    Любой полный заострённый веер размерности 2 политопален.
	}
	\end{proposition}
	
	\begin{definition}
	  \emph{Веером схождений} $\mathcal F_{v, \mathcal T}$ в вершине $v$ нормального 
		локально-конечного разбиения $\mathcal T$ называется набор конусов следующего вида:
		$$\mathcal F_{v, \mathcal T} = \left\{\cone_{v}(F)|~ v \in F, F \hbox{ --- грань разбиения} 
		  \mathcal T\right\}$$
	\end{definition}
	
	Легко видеть, что для произвольной грани $F \subset \mathcal T$, содержащей вершину $v$, 
	конус $\cone_{v}(F)$ задаётся в виде $\{v + \sum^k_{i=1}\alpha_i\overrightarrow{e_i}, 
	  \alpha_i \geqslant 0\}$, где $\overrightarrow{e_i}$ --- направленные рёбра грани $F$, выходящие
	из вершины $v$.
	
	Отсюда легко вывести, что веер схождений в вершине является полным заострённым $d$-мерным веером.
		
	\begin{definition}\label{definitionTightFansMeetingFan}
	  Пусть задано локально-конечное разбиение $\mathcal T$ и $n$-мерная грань $V \subset \mathcal T$.
		
		Через внутреннюю точку $v \in V$ проведём трансверсальное к $\aff{V}$ подпространство $S$ 
	  (размерности $d - n$). Тогда $S$ пересекает $V$ по единственной точке $v$, а любую $k$-грань, 
	  содержащую $V$ --- по выпуклому	многограннику размерности $k-n$ с вершиной $v$. 
		
		Таким образом, пересечение $S$ и $\mathcal T$ является нормальным $(d-n)$-мерным разбиением 
		$\widehat{\mathcal T}$ подпространства $S$ с вершиной $v$. 
		
		\emph{Веером схождений} $\mathcal F_{V, \mathcal T}$ в грани $V$ разбиения $\mathcal T$ 
		называется веер схождений в вершине $v$ разбиения $\widehat{\mathcal T}$.
	\end{definition}
	
	\begin{definition}\label{definitionVoronoiLinearPart}
	  \emph{Линейной частью} аффинного подпространства $L$ в $\mathbb E^d$ называется линейное 
		подпространство $L - \mathbf{p}$, где $\mathbf{p}$ --- произвольная точка из $L$.
		
		Линейная часть подпространства $L$ обозначается $\lin{L}$.
	\end{definition}
	
	Очевидно, что такой веер зависит только от линейной части $\lin{S}$ подпространства $S$ и не 
	зависит от выбора внутренней точки $v$ в грани $V$.
	
	Кроме того, если два веера схождений в $V$ построены по двум различным подпространствам $S_1$ и 
	$S_2$, то их можно перевести друг в друга аффинным преобразованием, сохраняющим аффинную 
  оболочку $\aff V$. То есть определение веера схождений в грани корректно и однозначно 
	с точностью до аффинного преобразования.

\subsection{Нормальные веера и полярные многогранники.}
	
	\begin{definition}\label{definitionTightFansPolar}
	  Выберем произвольную точку $p$ внутри выпуклого $d$-мерного многогранника $P \subset \mathbb 
		  E^d$. \emph{Полярным многогранником} для $P$ называется следующее точечное 
		множество:
		$$\polar{P}_p ~:= ~\{c \in \mathbb R^d: (c-p)\cdot(x-p) \leqslant 1 \hbox{ для всех } x\in P\}$$
		
		Точку $p$ будем называть \emph{центром полярности} многогранника $\polar{P}$.
	\end{definition}
	
  Хорошо известно, что $\polar{P}$ является выпуклым $d$-мерным многогранником и также содержит 
	точку $p$ внутри. Кроме того, известны важные свойства \cite[Глава 2]{Ziegler}:
    
  \begin{proposition}\label{propositionBasicProp}~{
    \begin{enumerate}
      \item $\dubpolar{P} = P$ 
			\item $k$-граням многогранника $P$ ($0\leqslant k \leqslant d-1$) однозначно соответствуют 
			  $(d-1 - k)$-грани многогранника $\polar{P}$. Для грани $F$ соответствующая ей грань в 
				$\polar{P}$ обозначается как $\polface{F}$.
			\item $\dubpolface{F} = F$
      \item Для двух граней $F, G \subset P$ грань $F$ является гранью $G$ тогда и только тогда, 
			  когда выполнено обратное включение: грань $\polface{F}$ содержит грань $\polface{G}$. 
    \end{enumerate}
  }
  \end{proposition}
	
	Последнее свойство предложения \ref{propositionBasicProp} означает, что полярное соответствие
	граней $\diamond$ сохраняет инцидентности граней многогранника, но меняет порядок их включения.
	
	\begin{definition}
	  Пусть задана точка $p \in \mathbb E^d$ и не содержащая $p$ гиперплоскость $\mathcal H \subset
		  \mathbb E^d$. \emph{Полюсом} гиперплоскости $\mathcal H$ относительно центра $p$ называется 
		такая точка $v_{\mathcal H} = v_{p, \mathcal H} \in \mathbb E^d$, что для любой точки $t \in 
		  \mathcal H$ выполнено $(v_{\mathcal H} - p)\cdot (t - p) = 1$.
	\end{definition}
	
	Легко видеть, что полюс $v_{p, \mathcal H}$ определён однозначно, и что $v_{p, \mathcal H}$ лежит 
	на луче, ортогональном гиперплоскости $\mathcal H$ и пересекающем её.
	
	В \cite[Глава 2]{Ziegler} доказано следующее предложение:
	
	\begin{proposition}\label{propositionTightFansCorrespondence}~{
	  Пусть задан $d$-мерный многогранник $P \subset \mathbb E^d$ и его полярный многогранник 
		$\polar{P}_p$. Для произвольной гиперграни $H \subset P$ соответствующая ей вершина 
		$\polface{H} \in \polar{P}_p$ совпадает с полюсом $v_{p, \aff{H}}$ аффинной оболочки 
		гиперграни $H$.
	}
	\end{proposition}
	
	Два полярных многогранника $\polar{P}_1$ и $\polar{P}_2$, построенные относительно разных центров 
	$p_1$ и $p_2$, в общем случае не являются аффинно-эквивалентными. Однако для каждого 
	$\polar{P}_i$ инцидентности его граней однозначно определены инцидентностями граней в $P$.
	Значит	$\polar{P}_1$ и $\polar{P}_2$ имеют одинаковую комбинаторную структуру, то есть 
	\emph{комбинаторно эквивалентны}.
	
	\begin{definition}
	  Пусть задан выпуклый $d$-мерный многогранник $P \subset \mathbb E^d$. \emph{Нормальным веером}
	  $\normalfan{P}$ к $P$ называется веер граней его полярного многогранника $\polar{P}_p$
		для произвольного центра полярности $p \in \int{P}$.
	\end{definition}
	
	По данному определению, нормальный веер зависит от выбора центра полярности. Следующее утверждение 
	является фольклорным в теории вееров. Его доказательство приводится здесь для полноты изложения. 
	
	\begin{lemma}\label{lemmaTightFansNormalFanUnique}~{
	  Нормальный веер $\normalfan{P}$ однозначно (с точностью до параллельного переноса) задаётся
		выпуклым многогранником $P$.
	}
	\end{lemma}
	
	\smallskip
	\noindent {\it Доказательство.} 
	  Полный заострённый веер однозначно задаётся своими одномерными рёбрами и инцидентностями 
		граней. Для двух различных полюсов $p_1, p_2 \in \int{P}$, одномерные рёбра вееров
		$\mathcal N_1(P)$ и $\mathcal N_2(P)$ совмещаются параллельным переносом на вектор 
		$\overrightarrow{p_1p_2}$. Инцидентности граней этих вееров такие же, как у $\polar{P}_1$ и 
		$\polar{P}_2$ соответственно, а у них инцидентности одинаковые.
	\hfill $\Box$ \par\bigskip
	
	Другим стандартным определением нормального веера является определение через линейные функционалы
	в $\left(\mathbb E^d\right)^*$, которое можно найти в \cite[Пример 7.3]{Ziegler}. 

	Далее будем обозначать через $N_F$ конус нормального веера $\normalfan{P}$, который соответствует 
	грани $F \subseteq P$. Отметим, что если $\normalfan{P}$ задан через центр полярности $p \in 
	\int{P}$, то $N_F$, по определению, совпадает с конусом $\cone_p\left(\polface{F}\right)$.

	\begin{lemma}\label{lemmaTightFansOrthogSubspace}~{\\
	  а)\quad Пусть $F$ --- произвольная $k$-грань $P$ ($k \geqslant 1$), тогда $F$ ортогональна
		$N_F$ --- соответствующей грани нормального веера.\\
		б)\quad Пусть $\polar{P}$ --- полярный многогранник, построенный по полюсу $p \in \int{P}$, 
		тогда $\cone_p{F}$ ортогонален соответствующей грани $\polface{F}$ полярного многогранника.
	}
	\end{lemma}
	
	\smallskip
	\noindent {\it Доказательство.} 
	  Докажем сначала пункт а) леммы. Построим $\normalfan{P}$ как веер граней $\polar{P}$ для 
		произвольного полюса $p \in \int{P}$. Обозначим $H_1, \ldots, H_n$ все гиперграни $P$, 
		содержащие $F$. Тогда по предложению \ref{propositionTightFansCorrespondence} имеем 
		$\overrightarrow{p\polface{H_i}} \bot H_i \supset F$, откуда 
		$\overrightarrow{p\polface{H_i}} \bot \aff{F}$. Остаётся заметить, что
	  $$N_F = \cone_p\left(\polface{F}\right) = \cone_p\left(\overrightarrow{p\polface{H_1}}, \ldots, 
		  \overrightarrow{p\polface{H_n}}\right) \bot \aff{F}.$$
		
		Пункта б) следует из соотношений $\dubpolar{P} = P, \dubpolface{F} = F$ (предложение 
		\ref{propositionBasicProp}) и того, что $\mathcal F(P)$ --- веер граней для $P$ --- совпадает с 
		$\normalfan{\polar{P}} := \mathcal F\left(\dubpolar{P}\right)$.
	\hfill $\Box$ \par\bigskip

\section{Тесные веера}
\subsection{Тесные веера, их редукция и разложимость.}
  \begin{definition}
	  Назовём нормальное разбиение $\mathcal T$ пространства $\mathbb E^d$ \emph{локально 
		симметричным} в грани $F$, если найдётся шаровая окрестность некоторой внутренней точки $p \in 
		  F$, внутри которой данное разбиение центрально-симметрично относительно $p$.
	\end{definition}

  \begin{lemma}\label{lemmaLocalSymm}~{
    Локальная симметрия нормального разбиения пространства $\mathbb E^d$ в грани $F$ определено 
		корректно и не зависит от выбора внутренней точки грани.
  }
  \end{lemma}
  
  \smallskip
	\noindent {\it Доказательство.} Рассмотрим такие шаровые окрестности двух различных внутренних
	  точек $p_1, p_2$ из $F$, что радиусы этих окрестностей равны, а сами окрестности не пересекают 
		других граней разбиения, кроме содержащих $F$. Такие окрестности, вместе с частью разбиения, 
		попадающей внутрь, совмещаются параллельным переносом. Значит они центрально-симметричны или 
		не центрально-симметричны одновременно.
	\hfill $\Box$ \par\bigskip
	
	\begin{definition}
    Грань $f$ полного веера будем называть \emph{стандартной}, если данный веер в ней локально 
		симметричен.
	\end{definition}
	
	\begin{definition}
	  Два конуса $c_1$ и $c_2$ (не обязательно полной размерности) в $\mathbb E^d$ называются
		\emph{локально симметричными} в общей грани $f$, если найдётся шаровая окрестность 
	  некоторой внутренней точки $p \in f$, внутри которой данные два конуса центрально симметричны 
		друг другу относительно $p$. 
	\end{definition}
	
	\begin{definition}
	  \emph{Тесным веером} называется такой полный заострённый $d$-мерный веер в $\mathbb E^d$, что 
		пересечение произвольных двух его $d$-мерных конусов $C_i, C_j$ является стандартной гранью 
		веера, в которой эти конусы локально симметричны.
	\end{definition}
	
	\begin{lemma}\label{lemmaTightFansSymmetry}~{
	  Пусть задан $d$-мерный полный веер $\mathcal F$ и $f$ --- некоторая его грань. Пусть грани 
		$f_1$ и $f_2$ этого веера содержат $f$ и им соответствуют грани $\widehat{f_1}$ и 
		$\widehat{f_2}$ веера схождений $\mathcal F_f$ в грани $f$. Тогда $f_1$ и $f_2$ локально 
		симметричны тогда и только тогда, когда $\widehat{f_1}$ и $\widehat{f_2}$ локально 
		симметричны.
	}
	\end{lemma}
	
	\smallskip
	\noindent {\it Доказательство.}
	  Рассмотрим веер схождений $\mathcal F_f$ как сечение веера $\mathcal F$ подпространством $L$ 
		трансверсальным к грани $f$ и проходящим через её внутреннюю точку (см. определение 
		\ref{definitionTightFansMeetingFan}).
		
		Пусть $f_1 \bigcap f_2 = g \subset \mathcal F$. Тогда $f \subset g$ и $\widehat{f_i} = f_i \bigcap L$, 
		$\widehat{g} = g \bigcap L$.
		
		Рассмотрим произвольную внутреннюю точку $p \in \int{\widehat{g}} \subset \int{g}$. Выберем 
		шаровую окрестность $\mathcal B_{\varepsilon}(p)$, которая не пересекает других граней из
		$\mathcal F$, кроме содержащих $g$. Отсюда следует, что $\mathcal B_{\varepsilon}(p)$
		пересекает $g$ только по внутренним точкам. Тогда шаровая окрестность $\widehat{\mathcal 
		  B_{\varepsilon}}(p)$ в $\mathcal F_f \subset L$ является центральным сечением 
		$\mathcal B_{\varepsilon}(p)$ подпространством $L$.
		
		Если $f_1$ и $f_2$ симметричны внутри $\mathcal B_{\varepsilon}(p)$ для некоторого $\varepsilon 
		  > 0$, то $\widehat{f_1} = f_1 \bigcap L$ и $\widehat{f_2} = f_2 \bigcap L$ симметричны в
		$\widehat{\mathcal B_{\varepsilon}}(p) = \mathcal B_{\varepsilon}(p) \bigcap L$. В одну сторону
		лемма доказана.
		
		Обратно, пусть $\widehat{f_1}$ и $\widehat{f_2}$ симметричны в $\widehat{\mathcal 
		  B_{\varepsilon}}(p)$ относительно точки $p$ для некоторого $\varepsilon > 0$.
			
		По выбору, произвольная грань $h \subset \mathcal F$, пересекающая $\mathcal B_{\varepsilon}(p)$,
		содержит грань $g$, а значит и грань $f$.  В частности, все гиперграни из $\mathcal F$, 
		пересекающие шар $\mathcal B_{\varepsilon}(p)$, содержат его центр $p$. Отсюда следует, что
		$h \bigcap \mathcal B_{\varepsilon}(p) = \cone_p{h} \bigcap\mathcal B_{\varepsilon}(p)$.
		
		Кроме того, $\aff{f} \subset \aff{g} \subset \aff{h}$ и, следовательно, $\lin{f} \subset 
		  \lin{g} \subset \lin{h}$. Отсюда легко видеть, что $\cone_p{h} = \cone_p{\widehat{h}}\oplus
			\lin{f}$. Следовательно $h \bigcap \mathcal B_{\varepsilon}(p) = \left(\cone_p{\widehat{h}}
			\oplus\lin{f}\right) \bigcap \mathcal B_{\varepsilon}(p)$. 
		
		Условие, что $\widehat{f_1}$ и $\widehat{f_2}$ симметричны в $\widehat{\mathcal 
		  B_{\varepsilon}}(p)$ относительно $p$ означает, что $\widehat{f_1} \bigcap \widehat{\mathcal 
		  B_{\varepsilon}}(p)$ симметрично $\widehat{f_2} \bigcap \widehat{\mathcal B_{\varepsilon}}(p)$
		в точке $p$. Значит, $\cone_v{\widehat{f_1}}$ симметрично $\cone_v{\widehat{f_2}}$ в точке $p$.
			
		Получаем, что  $f_1 \bigcap \mathcal B_{\varepsilon}(p) = \left(\cone_p{\widehat{f_1}}\oplus
		  \lin{f}\right) \bigcap \mathcal B_{\varepsilon}(p)$ симметрично
		$\left(\cone_p{\widehat{f_2}}\oplus\lin{f}\right) \bigcap \mathcal B_{\varepsilon}(p) = 
		  f_2 \bigcap \mathcal B_{\varepsilon}(p)$ в точке $p$. Что и 
		требовалось.
	\hfill $\Box$ \par\bigskip
	
	Следующие две теоремы задают правила сведения тесных вееров к тесным веерам меньшей размерности.
	
	\begin{theorem}\label{theoremTightFansSection}~{
	  Пусть задан тесный веер $\mathcal F$ и $f$ --- некоторая его грань, тогда веер схождений 
	  $\mathcal F_f$ в грани $f$ также является тесным.
	}
	\end{theorem}
	
	\smallskip
	\noindent {\it Доказательство.}
	  Пусть $\widehat{C}_i, \widehat{C}_j$ --- произвольные две ячейки полной размерности в $\mathcal
		  F_f$, а $\widehat{H}_{ij} = \widehat{C}_i \bigcap \widehat{C}_j$. По определению веера 
		схождений, $\widehat{H}_{ij}$ является сечением грани $H_{ij} = C_i \bigcap C_j$ 
		подпространством $L \supset \mathcal F_f$.
		
		Очевидно, $f$ принадлежит $H_{ij}$. По определению тесного веера, $H_{ij}$ --- стандартная 
		грань, и $C_i$ локально симметрично $C_j$ в $H_{ij}$. Из леммы \ref{lemmaTightFansSymmetry}
		следует, что $\widehat{C}_i$ локально симметрично $\widehat{C}_j$ в $\widehat{H}_{ij}$.
		
		Кроме того, из определения стандартной грани следует, что полномерные ячейки веера $\mathcal F$
		распадаются на пары симметричных в $H_{ij}$. Из леммы \ref{lemmaTightFansSymmetry} следует, 
		что полномерные ячейки веера $\mathcal F_f$ также разобьются на пары симметричных в 
		$\widehat{H}_{ij}$.
		
		Значит, пересечение произвольных двух ячеек $\widehat{C}_i, \widehat{C}_j$ полной размерности 
		в $\mathcal F_f$ является стандартной гранью в $\mathcal F_f$, в которой данные ячейки 
		локально симметричны. По определению, $\mathcal F_f$ --- тесный веер.
	\hfill $\Box$ \par\bigskip
	
	\begin{definition}\label{definitionVoronoiMinkovskiSum}
    Напомним, что \emph{суммой Минковского} множеств\\ $Q_1, Q_2, \ldots, Q_k \subset \mathbb E^n$ 
	  называется множество точек 
		$$Q = \sum^k_{i = 1}{Q_i} := \left\{\ib{x}\: \bigl |\: \ib{x} = 
	    \sum^k_{i = 1}{\ib{q}_i}, \ib{q}_i \in Q_i\right\}$$ 
		Если слагаемые расположены в трансверсальных подпространствах, то такая сумма называется 
		\emph{прямой суммой Минковского}.
  \end{definition}
	
	\begin{definition}
	  \emph{Прямой суммой вееров} (см. \cite{Ziegler}) $\mathcal F_1 \subset \mathbb E^{k_1}$ и 
		$\mathcal F_2 \subset \mathbb E^{k_2}$ называется следующий набор конусов в 
		$\mathbb E^{k_1}\times \mathbb E^{k_2}$:
  	$$\mathcal F_1\oplus \mathcal F_2 := \left\{C_1\oplus C_2 \subset \mathbb E^{k_1}\times 
		  \mathbb E^{k_2}: C_1 \in \mathcal F_1, C_2\in \mathcal F_2\right\},$$
	  где $C_1, C_2$ --- грани соответствующих вееров, а $\oplus$ --- прямая сумма Минковского 
		подмножеств евклидова пространства.		
	\end{definition}
	
	Легко видеть, что прямая сумма полных заострённых вееров также будет полным заострённым веером.
	
  \begin{theorem}\label{theoremTightFans}~{
	  Пусть задан тесный веер $\mathcal F$, который представляется в виде прямой суммы вееров
	  $\mathcal F_1\oplus \mathcal F_2$. Тогда $\mathcal F_1$ и $\mathcal F_2$ --- полные, заостроённые
	  и тесные веера.
	}
	\end{theorem}
	
	\smallskip
	\noindent {\it Доказательство.} 
	  Веера $\mathcal F_1$ и $\mathcal F_2$ должны быть полными, иначе прямая сумма не совпадает с 
		$\mathbb R^d = \mathbb R^{k_1} \times \mathbb R^{k_2}$. Кроме того, грани 
	  суммарного веера равны сумме граней, значит вершина могла получиться только как сумма нульмерных
	  вершин. Значит, слагаемые прямой суммы --- заострённые веера. 
	  
	  Рассмотрим произвольную $k_2$-мерную ячейку $C_2$ веера $\mathcal F_2$. Конусы множества 
		$\mathcal L = \left\{C_1\oplus C_2|~ C_1 \hbox{ --- грань в } \mathcal F_1\right\}$ суть все 
		конусы из веера $\mathcal F_1\oplus \mathcal F_2$, которые содержат $k_2$-грань $v_1 \oplus 
		  C_2$ (тут $v_1$ --- вершина $\mathcal F_1$). 
			
		Произвольное трансверсальное сечение к этой грани, проведённое через её внутреннюю точку 
		$v_1 \oplus p, p \in \int{C_2}$, локально является веером схождений в $v_1 \oplus C_2$ и 
		аффинно-эквивалентно вееру $\mathcal F_1$. 
		
		С другой стороны, по теореме \ref{theoremTightFansSection}, веер схождений в грани $v_1 \oplus 
		  C_2$ тесного веера $\mathcal F_1\oplus \mathcal F_2$ также является тесным. Значит и веер 
		$\mathcal F_1$ тесный, так как локальные симметрии в гранях сохраняются при аффинных 
		преобразованиях. Аналогично второй веер $\mathcal F_2$ также тесный.
	\hfill $\Box$ \par\bigskip

  Далее мы даём классификацию тесных вееров размерности 2.

  \begin{lemma}\label{lemma2dTight}~{
    Существует 2 различных комбинаторных типа тесных вееров в размерности 2: разбиение плоскости, 
    получаемое проведением трёх лучей с общей вершиной (и каждым из полученных углов меньше 
		развёрнутого); разбиение плоскости, получаемое проведением двух пересекающихся прямых.
  }
  \end{lemma}
  
  \smallskip
	\noindent {\it Доказательство.}
	  Ячейки полной размерности полных заострённых вееров размерности 2 --- это плоские углы меньшие
		развёрнутого. Значит, плоский веер содержит не менее трёх углов. 
		
		Плоский веер с тремя углами (меньшими развёрнутого) является тесным, так как пересечение любых 
		двух из этих углов --- луч. Как легко проверить, два плоских угла веера всегда локально 
		симметричны в луче пересечения.
	  
	  Иначе, в веере $k \geqslant 4$ угла. Обозначим эти углы $A_1, A_2, \ldots, A_k$, занумерованные
		против часовой стрелки. Тогда $A_i$ для всякого $0 \leqslant i \leqslant k$, пересекается с 
		$A_{i-1}$ и $A_{i+1}$ по лучу, а остальными углами --- лишь по общей вершине. Тогда, так как
		веер тесный, то угол $A_i$ локально симметричен каждому из углов $A_{i+2}, A_{i+3}, 
	    \ldots, A_{i-2}$ в общей вершине. Что, очевидно, невозможно, если от $A_{i+2}$ до $A_{i-2}$ 
		углов более одного. То есть при $k > 4$. 
		
		При $k = 4$ получаем, что $A_1$ локально-симметричен $A_3$ в общей вершине. Значит, эти углы 
		являются вертикальными и образованы пересечением двух прямых.
	\hfill $\Box$ \par\bigskip

\subsection{Комбинаторная классификация тесных вееров размерности 3.}
  В разделе мы приводится классификация трёхмерных тесных вееров. Показано, что существует
  всего 5 различных комбинаторных типов тесных вееров. Они в точности совпадают с пятью типами 
	Делоне \cite{Delone} --- типами схождений параллелоэдров в гранях размерности $(d-3)$. 
	
	Теорему \ref{theoremTightFansClassification} другим методом доказал А. Магазинов в 
	\cite{Magazinov3Tight}\footnote{Доказательства были получены 
	приблизительно одновременно при совместных исследованиях авторов в Queen's University, 
	г. Кингстон, Канада.}
	
	\begin{theorem}\label{theoremTightFansClassification}~{
	  Существует 5 различных комбинаторних типов трёхмерных тесных вееров. Комбинаторно они совпадают
		с веерами граней трёхмерных многогранников: октаэдра, четырёхгранной пирамиды, куба, тетраэдра 
		и трёхгранной бипирамиды. 
	}
	\end{theorem}
	
	Следует отметить, что приведённая классификация является комбинаторной, то есть верно отражает
	инцидентности всех граней, но не гарантирует, что конкретный веер можно получить именно как
	как конус над поверхностью реального многогранника. Этот факт верен, но будет доказан позже. 
	Для его доказательства потребуется, в том числе, и эта теорема.
  
  Рассмотрим циклическую последовательность трёхмерных (то есть полномерных) попарно различных 
  конусов $[C_1, C_2, \ldots, C_k]$ из данного веера $\mathcal F$, в которой любые два 
  последовательных конуса смежны по грани (размерности 2). Такую цепочку конусов можно 
  рассмотреть как конус над двумерным кольцом. Действительно: рассмотрим двумерную сферу 
  $S = S^2_v$ с центром в вершине веера $v$. Произвольному выпуклому конусу $C$ с вершиной $v$ 
  однозначно сопоставим выпуклый сферический многогранник $\widetilde{C} = C \bigcap S$ на 
  сфере. Тогда цепочке конусов однозначно соответствует цепочка сферических многоугольников
  $[\widetilde{C_1}, \widetilde{C_2}, \ldots, \widetilde{C_k}]$. По этой цепочке построим 
  (сферическую) замкнутую ломаную $L = [p_1, p_2, \ldots, p_k]$:  на одномерном ребре между 
  $\widetilde{C_i}$ и $\widetilde{C_{i+1}}$ выберем произвольную внутреннюю точку $p_i$ и 
  соединим их (сферическими) отрезками (то есть дугами центральных сечений сферы). Из 
  выпуклости многоугольников $\widetilde{C_i}$ следует, что ломаная $L$ лежит строго внутри 
  объединения ячеек $\bigcup\limits^k_{i = 1}\widetilde{C_i}$. Кроме того, так как все 
  конусы цепочки различны, значит все многоугольники сферической цепочки также различны. Отсюда 
  вытекает, что $L$ --- простая кривая без самопересечений. $L$ делит поверхность сферы $S$ на 
  три части --- точки кривой $l$ и два двумерных дополнительных множества. Каждый конуса $C_j$, 
  не принадлежащий цепочке, пересекает ровно одно из этих двух сферических множеств (и не 
  пересекает ломаную на ней). Одно из этих двух множеств, условно, назовём \emph{северным} 
  множеством конусов, а второе --- \emph{южным}. Очевидно, южное и северное множества не зависят от 
  конкретного выбора точек ломаной $[p_1, \ldots, p_k]$. Кроме того, легко привести пример, когда 
  хотя бы одно из этих двух множеств конусов (северного и южного) было пустым. Если оба эти
  множества непусты, то такую замкнутую цепочку конусов $[C_1, \ldots, C_k]$ будем называть 
  \emph{разделющей}.
  
  \begin{lemma}\label{lemmaTightFansSeparating}~{
    Если цепочка $C_1, \ldots, C_k$ в тесном веере является разделяющей, то вершина веера 
    является его центром симметрии (то есть стандартна), цепочка симметрична себе в этой вершине, 
    а дополнение к цепочке состоит ровно из двух симметричных друг другу (также в вершине) конусов.
  }
  \end{lemma}
  
  \smallskip
	\noindent {\it Доказательство.} 
	  Выберем $d$-конус $U_1$ не из цепочки и произвольный конус $U_2$ из второго дополнительного
	  множества. Тогда эти два конуса имеют лишь общую вершину $v$. Значит она стандартна, а 
	  конусы $U_1$ и $U_2$ в ней симметричны. Если второе дополнительное множество содержит ещё
	  хотя бы один конус $U_3$, то $U_1$ и $U_3$, исходя из тех же рассуждений, также симметричны
	  в вершине $v$, что невозможно, так как $U_2$ и $U_3$ различны по выбору. Аналогично, в первом
	  дополнительном множестве, не может быть других конусов кроме $U_1$. Множество конусов цепочки,
	  таким образом, состоит и множества конусов веера без подмножества $\{U_1, U_2\}$. Оба эти 
	  множества конусов симметричны в вершине, значит и цепочка также симметрична.
	\hfill $\Box$ \par\bigskip
	
	\begin{corollary}\label{corollaryTightFansNo3}~{
	  Тесный веер не содержит разделяющего цикла длины 3.
	}
	\end{corollary}
	
  Аналогично определяется \emph{частично разделяющая} цепочка конусов. Если в цепочке полномерных 
  конусов $[C_1, \ldots, C_k]$ данного веера $\mathcal F$ каждая пара $C_i, C_{i+1}$ (для $i = 1,
    2, \ldots, (k-1)$) смежна по двумерной грани, а пара $C_k, C_1$ смежна по одномерному ребру 
  $l$, то для такой цепочки аналогично строится сферическая ломаная $L$, вершины $p_1, \ldots,
    p_{k-1}$ которой выбираются как внутренние точки соответствующих сферических рёбер 
  $\widetilde{C_i} \bigcap \widetilde{C_{i+1}}$, а последняя вершина $p_k$ --- пересечение сферы
  с ребром $l$, то есть с одномерной гранью $\widetilde{C_k} \bigcap \widetilde{C_1}$. Относительно
  этой ломаной также строятся северное и южное множества конусов $C_j$ (из тех, что не принадлежат
  цепочке). Отличие от рассмотренной ситуации в том, что эти конусы могут иметь пересечение с 
  данной ломаной. Но лишь по единственной точке --- по $p_k$. Если оба множества конусов --- 
  северное и южное, для данной цепочки непусты, то мы будем называть её \emph{частично 
  разделяющей}. Легко доказывается следующее предложение.
	
	\begin{proposition}\label{propositionTightFansSemiSeparating}~{
    Если цепочка $C_1, \ldots, C_k$ в тесном веере является частично разделяющей ($C_k \bigcap C_1$ 
    --- одномерное ребро $l$), а конусы $C_i$ и $C_j$ выбраны одно из северного, другой --- из 
    южного относительно этой цепочки множеств, то $C_i$ и $C_j$ пересекаются либо по вершине, либо 
    по ребру $l$.
  }
  \end{proposition}

  \begin{lemma}\label{lemmaTightFansCubeType}~{
    Если вершина трёхмерного тесного веера стандартна, а сам веер содержит хотя бы одно стандартное 
    (одномерное) ребро, то все его рёбра стандартны. Такой веер представляется в виде прямой суммы 
    трёх одномерных тесных вееров или, иначе, двумерного тесного веера со четырмя рёбрами и 
		единственного одномерного тесного веера. Комбинаторно этот веер совпадает с веером граней 
		октаэдра.
  }
  \end{lemma}
  
	\smallskip
	\noindent {\it Доказательство.} 
	  Обозначим стандартное ребро как $l$, вершину как $v$, а веер как $\mathcal F$. По лемме 
	  \ref{lemmaTightFansSymmetry} веер схождений $\mathcal F_l$ в $e$ стандартен. Для трёхмерного 
	  $\mathcal F$ веер схождений в ребре $l$ будет двумерен. Более того, стандартное ребро $l$ 
	  перейдёт в стандартную вершину. По лемме \ref{lemma2dTight} существует лишь один тип двумерного 
	  тесного веера, в котором вершина стандартна. Для этого типа в вершине сходятся 4 двумерные 
	  ячейки. Обозначим их $\widehat{C_1}, \widehat{C_2}, \widehat{C_3}, \widehat{C_4}$ согласно 
	  обходу (любому из двух) вокруг вершины веера $\mathcal F_e$. Они соответствуют 3-ячейкам 
	  $C_1, C_2, C_3, C_4$, сходящимся в $e$.
	  
	  Пусть $C^*_i$ --- ячейка в $\mathcal F$ центрально-симметричная $C_i$ в вершине $v$ (стандартна
	  по условию), $l^*$ --- ребро симметричное $l$ в $v$. Тогда $C*_1 \bigcap C_2$ --- стандартная 
	  грань веера, не совпадающая с $v$. Покажем, что эта грань также не может быть двумерной. Из
	  локальной симметрии в $l$ следует, что двумерные грани $C_2\bigcap C_3$ и $C_4\bigcap C_1$ 
	  лежат в одной двумерной плоскости $\alpha$, содержащей грань $e$, а также $v \in e$. При 
	  локальной симметрии в вершине $v$ грань $C_4\bigcap C_1$ перейдёт (локально) в грань 
	  $C^*_4\bigcap C^*_1$ и продолжит находится в $\alpha$. Отсюда легко видеть, что ячейки $C^*_1$
	  и $C_2$ отделены двумерной плоскостью $\alpha$ (лежат по разные стороны от неё), и содержат в
	  ней свои двумерные грани $C^*_4\bigcap C^*_1$ и $C_2\bigcap C_3$ соответственно. Отсюда следует, 
	  что общая грань $C^*_1 \bigcap C_2$ совпадает с пересечением граней 
	  $\left(C^*_4\bigcap C^*_1\right) \bigcap \left(C_2\bigcap C_3\right)$. Эти грани, 
	  очевидно, различны. Если бы общая грань $C^*_1 \bigcap C_2$ была двумерна, то должна была
	  бы лежать также в $\alpha$ и, скажем, ячейчка $C_2$ не была бы строго выпуклой. 
	  
	  Таким образом, $C^*_1 \bigcap C_2$ --- это одномерная грань $l_2$, которая принадлежит $C_3$. 
	  Аналогично, показывается, что $C^*_1 \bigcap C_4$ --- это одномерная грань $l_4$, также 
	  принадлежащая ячейкам $C_3, C_4, C^*_2, C^*_1$. Отсюда следует, что двумерный угол 
	  $\langle l_2, l_4\rangle$, натянутый на одномерные конусы $l_2, l_4$ --- это пересечение
	  ячеек $C^*_1$ и $C_3$. Обозначим через $l^*_2$ и $l^*_4$ образы $l_2$ и $l_4$ при симметрии в
	  $v$. Тогда аналогично имеем $\langle l_2, l^*_4\rangle = C_2 \bigcap C^*_4$, а из симметрии
	  $\langle l^*_2, l^*_4\rangle = C_1 \bigcap C^*_3$ и $\langle l^*_2, l_4\rangle = C^*_2 \bigcap 
	    C_4$. Очевидно, все грани веера нами перечислены. Получаем, что данный веер $\mathcal F$ 
	  равен прямой сумме одномерных вееров $\{l, l^*\}, \{l_2, l^*_2\}$ и $\{l_4, l^*_4\}$.
	\hfill $\Box$ \par\bigskip

	Из определения следует, что произвольный трёхмерный конус в тесном веере --- пространственный 
	угол, ограниченный некоторым количеством плоских углов. Это количество не может быть менее 3, 
	иначе это противоречило бы строгой выпуклости конусов (так как веер заострённый и у каждого
	конуса есть вершина). Покажем, что это количество не может быть более 4.

	\begin{lemma}\label{lemmaTightFans4Neigbors}~{
	  Если произвольный полномерный конус $C$ в трёхмерном тесном веере имеет более трёх плоских 
	  граней, то их ровно 4, а конусы, смежные с $C$ по этим двумерным граням, либо все 
	  одновременно имеют общее ребро, либо никакие два из них не имеют общего ребра.
	}
	\end{lemma}
	
	\smallskip
	\noindent {\it Доказательство.} 
	   Обозначим двумерные грани $C$ через $H_1, H_2, \ldots, H_k$, в порядке их следования при 
	   обходе вокруг вершины. Рассмотрим произвольные две несоседние грани $H_i, H_j$ и 
	   соответствующие конусы $C_i, C_j$ данного веера такие, что $C_i\bigcap C = H_i, 
	     C_j\bigcap C = H_j$. $C_i$ и $C_j$ не могут пересекаться по грани размерности 2, иначе
	   $C, C_i, C_j$ --- разделяющий цикл длины 3, чего не бывает в тесном веере, согласно лемме
	   \ref{corollaryTightFansNo3}. То есть $C_i$ и $C_j$ пересекаются либо по ребру, либо по
	   вершине. 
	   
	   Пусть $C_i$ и $C_j$ пересекаются по некоторому ребру $l$, тогда $C, C_i, C_j$ --- частично 
	   разделяющий цикл. Рассмотрим произвольную пару $C_m, C_n$ --- конусов смежных с $C$ по 
	   двумерной грани, из которых один конус лежит в северном, другой --- в южном множестве 
	   относительно частично разделяющей цепочки $C, C_i, C_j$. По предложению 
	   \ref{propositionTightFansSemiSeparating} $C_m$ и $C_n$ пересекаются либо по вершине веера, 
	   либо по ребру $l$. Если они пересекаются по вершине, то по определению тесного веера, 
	   вершина веера стандартна. Так же как и ребро $l$. Но тогда по лемме 
	   \ref{lemmaTightFansCubeType} этот веер имеет тип схождения как в разбиении на параллелепипеды.
	   В частности, каждая ячейка имеет три плоские грани. Но по условию $C$ имеет не менее четырёх
	   плоских граней --- противоречие. Значит $C_m$ и $C_n$ пересекаются по ребру $l$. Мы доказали, 
	   что если среди соседей $C$ по плоским граням, есть хотя бы два (не подряд идущих), которые 
	   смежны по ребру, то все они имеют общее ребро. Иначе ни одной такой пары нет.
	\hfill $\Box$ \par\bigskip

	\begin{lemma}\label{lemmaTightFans4Types}~{
	  Существует ровно два различных комбинаторных типа трёхмерных тесных вееров, содержащих
		хотя бы одну ячейку, ограниченную ровно четырмя двумерными углами. Это, комбинаторно, веера 
		граней четырёхгранной пирамиды и куба.
	}
	\end{lemma}
	
	\smallskip
	\noindent {\it Доказательство.} 
	  Рассмотрим данный веер $\mathcal F$ и соответствующую 3-ячейку $C$, ограниченную 4 плоскими
		углами $H_1, \ldots, H_4$. По предыдущей лемме возможны два случая. Первый случай --- все 
		ячейки смежные с $C$ по плоским граням $H_i$ содержат общее ребро. Тогда, очевидно, других
		трёхмерных ячеек в данном веере нет и это в точности комбинаторный тип конуса над 
		четырёхгранной пирамидой.
		
		Во втором случае соответствующие соседние с $C$ ячейки $C_1, \ldots, C_4$ делятся на пары
		симметричных в вершине $v$, которая стандартна. Следовательно в данном веере нет стандартных
		одномерных рёбер (иначе по лемме \ref{lemmaTightFansCubeType} в данном веере нет ячеек с
		4 двумерными гранями). Таким образом, $l_i$ --- общее ребро двумерхных граней $H_i, H_{i+1}$,
		не стандартно. Согласно теореме \ref{theoremTightFansSection} веер схождений в $l_i$ также тесный, а 
		согласно классификации двумерных тесных вееров \ref{lemma2dTight}, при не стандартной грани 
		$l_i$, этот веер содержит 3 ячейки полной размерности (которые по построению являются сечениями
		ячеек, содержащих $l_i$). Значит $C, C_i, C_{i+1}$ --- это все три ячейки, сходящиеся в $l_i$.
		Очевидно отсюда, что $C_i$ пересекается с $C_{i+1}$ по двумерной грани. Для замкнутой цепочки
		ячеек $[C_1, C_2, C_3, C_4]$ северное и южное множество конусов будут непусты (одно из них 
		содержит $C$, другое --- $C^*$, симметричный $C$ в вершине). Значит эта цепочка разделяющая и 
		для неё выполняются условия леммы \ref{lemmaTightFansSeparating}. Значит кроме четырёх ячеек
		этой цепочки есть всего две другие --- уже указанные $C$ и $C^*$, а комбинаторный тип веера ---
		конус над поверхностью куба.
	\hfill $\Box$ \par\bigskip
	
	Мы готовы завершить доказательство основной теоремы данного раздела. Остался случай, когда 
	каждая трёхмерная ячейка конуса ограничена ровно тремя двумерными гранями. Если в веере нет
	ни одного стандартного ребра, то выберем произвольную ячейку $C$ и её грани обозначим $H_1, H_2, 
	  H_3$. По этим граням с $C$ смежны ячейки $C_1, C_2, C_3$. Тогда, как показывалось ранее,
	из нестандартности рёбер следует, что $C_1$ и $C_2$ смежны по грани, также как и пары ячеек
	$C_2, C_3$ и $C_3, C_1$. Если кроме перечисленных четырёх ячеек есть ещё какие-то, то $[C_1, 
	  C_2, C_3]$ --- разделяющая цепочка длины 3, чего не может быть по лемме 
  \ref{corollaryTightFansNo3}. Таким образом, в данном случае (нет стандартных рёбер) стандартный
  веер состоит из 4 ячеек, схождение имеет комбинаторный тип конуса над поверхностью тетраэдра.
	
	Наконец, пусть веер содержит стандартное ребро $l$. По свойству стандартного ребра коразмерности
	2, в нём сходятся 4 ячейки, попарно симметричные в этом ребре \ref{lemma2dTight}. Обозначим их 
	$C_1, \ldots, C_4$, согласно обходу вокруг этого ребра. У произовольной из этих четырёх ячеек, 
	скажем, у $C_2$, всего 3 двумерные грани, две из которых она делит с $C_1$ и $C_3$, а ещё одна
	грань общая с некоторой ячейкой не из данных четырёх. Обозначим через $l_1$ второе ребро 
	пересечения $C_1$ с $C_2$ (которое вместе с $l$ задаёт грань $C_1 \bigcap C_2$). Аналогично 
	$C_2$ и $C_3$ задают ребро $l_2$, $C_3, C_4$ --- $l_3$ и $C_4, C_1$ --- $l_4$. 
	
	Среди рёбер $l_1, \ldots, l_4$ есть стандартное. Действительно, иначе ячейка $C^*$, смежная с 
	$C_1$ по ``третьей'' грани $\langle l_4, l_1\rangle$, также смежна по двумерной грани с $C_2$ 
	(так как $l_1$ не стандартно), то есть $C^*$ содержит грань $\langle l_1, l_2\rangle$, в 
	частности --- содержит ребро $l_2$. Аналогично получим, что $C^*$ содержит все рёбра $l_i$ и, 
	очевидно, имеет 4 двумерные грани. Что неверно по предположению. То есть, не умаляя общности, 
	ребро $l_4$ стандартно. В нём сходятся ячейки $C_1, C_4, C^*_1, C^*_4$ ($C_i$ и $C^*_i$ 
	локально в нём симметричны). Из локальной симметрии в рёбрах следует, что грани $C^*_1\bigcap 
	  C^*_4$, $C_4\bigcap C_1$ и $C_3\bigcap C_2$ лежат в одной плоскости $\alpha$. Причём ячейки 
	$C^*_1, C_4, C_3$ --- в одном полупространстве относительно $\alpha$, а $C^*_4, C_1, C_2$ --- 
	в другом. Тогда ячейки $C^*_1$ и $C_2$ отделены друг от друга плоскостью $\alpha$, имеют в 
	этой плоскости двумерные (различные) грани. Значит их пересечение лежит в $\alpha$, но 
	не может быть гранью размерности 2. Также как по предположению не может быть вершиной (иначе 
	она была бы стандартна). Значит это ребро, причём ребро $l_2$. Отсюда следует, что $l_2$ 
	стандартно. Получили, что $\alpha$ разбивается на грани (двумерные углы) $C^*_1\bigcap C^*_4$, 
	$C_4\bigcap C_1$ и $C_3\bigcap C_2$. Если ячейки $C^*_1, C_4, C_3$ не пересекаются в общем 
	ребре, то образуют разделющий цикл длины 3, что невозможно. Значит, пересекаются в ребре $l_3
	  \in C_3 \bigcap c_4$. Аналогично ячеки $C^*_4, C_1, C_2$ пересекаются в общем ребре $l_1$.
	Тем самым мы получили последний пятый комбинаторный тип трёхмерных тесных вееров --- тип конуса
	над поверхностью бипирамиды. Отметим, что $l_1 \in \langle l_1, l_4 \rangle \bigcap \langle l_1, l_2
	  \rangle$, то есть лежит в одномерном пересечении плоскостей $\aff{l_1, l_4}$ и $\aff{l_1, 
		l_2}$. Из локальной симметрии следует, что это те же плоскости, что  $\aff{l_3, l_4}$ и 
  $\aff{l_3, l_2}$. То есть $l_1$ и $l_3$ лежат на одной прямой. Тем самым доказано предложение
	
	\begin{proposition}\label{propositionTightFans4Types}~{
		В трёхмерном тесном веере комбинаторного типа, соответствующего трёхгранной бипирамиде, 
		содержится всего 5 одномерных рёбер, три из которых стандартны и лежат в одной плоскости, а два 
		других лежат на одной прямой, трансверсальной этой плоскости. Геометрически этот веер равен 
		прямой сумме двумерного тесного веера с тремя рёбрами и единственного одномерного тесного 
		веера.
	}
	\end{proposition}
	
	Теорема о комбинаторной классификации трёхмерных тесных вееров доказана.

  \section{Веера и канонические нормировки}
	\subsection{Политопальность трёхмерных тесных вееров.}
	
	Как мы отметили ранее, уже в размерности 3 существуют примеры неполитопальных вееров, то есть
	примеры полных заострённых вееров, которые невозможно представить в виде веера граней некоторого
	многогранника. Хотя мы и сформулировали классификацию трёхмерных тесных вееров в терминах 
	вееров граней, но эта классификация комбинаторная. То есть мы показали, что инцидентности граней
	этих вееров такие же, как у указанных вееров граней. Обязательно ли сами эти тесные веера 
	представляются (геометрически, а не комбинаторно), как веер граней над каким-то реальным
	многогранником в пространстве? Ответ на этот вопрос утвердительный.
	
	\begin{theorem}\label{theoremTightFansPolytopical}~{
	  Каждый трёхмерный тесный веер можно реализовать как веер граней некоторого многогранника.
	}
	\end{theorem}
	
	\smallskip
	\noindent {\it Доказательство.} 
	  В случае (комбинаторного) веера граней тетраэдра, достаточно на каждом луче отметить по точке
		$v_i$ и взять их выпуклую оболочку. Получится также тетраэдр, который будет (геометрически)
		порождать веер. 
		
		В случае трёхгранной бипирамиды веер содержит 3 стандартных одномерных ребра. Из стандартности
		следует, что они лежат в одной плоскости $\pi$. На каждом из них произвольно выберем по точке:
		$v_1, v_2, v_3$. Два оставшихся ребра, как отмечено в предложении 
		\ref{propositionTightFans4Types}, дополняют друг друга до прямой трансверсальной $\pi$. На них
		также выберем по точке: $v_4, v_5$. Очевидно, что $\conv{v_1, \ldots, v_5}$ --- бипирамида над
		$\conv{v_1, v_2, v_3}$. Случай комбинаторного веера граней октаэдра рассматривается аналогично, 
		с той разницей, что сначала мы выберем 4 стандартных ребра в одной двумерной плоскости.
		
		В случае четырёхугольной пирамиды одна из 3-ячеек веера ($C_1$) имеет 4 двумерные грани, а 
		остальные ячейки --- по три. $C_1$ --- строго выпуклый конус с вершиной. Тогда, во-первых,
		существует опорная плоскость $\pi$ к $C_1$, содержащая только вершину. Во-вторых, параллельное
		к $\pi$ сечение $C_1$ пересекает этот конус по четырёхугольнику $v_1v_2v_3v_4$. Выберем на 
		оставшемся ребре веера точку $v_5$ и $\conv{v_1, v_2, v_3, v_4, v_5}$ будет искомой
		порождающей пирамидой для данного веера.
		
		Оставшийся случай комбинаторного веера над кубом потребует более детального доказательства.
		Такой веер состоит из шести 3-ячеек, каждая из которых ограничена четырьмя двумерными гранями
		и, соответственно, четырьмя рёбрами. Вершина $v$ веера стандартна, то есть веер в ней локально
		симметричен. Для вершины это означает, что веер просто симметричен в $v$. Обозначим $r_1, r_2, 
		  r_3, r_4$ --- одномерные рёбра одной из шести ячеек (идущие в порядке обхода вокруг вершины).
		Пусть $l_1$ --- прямая пересечения плоскостей $\aff{r_1, r_2}$ и $\aff{r_3, r_4}$. $l_2$ ---
		прямая пересечения плоскостей $\aff{r_1, r_4}$ и $\aff{r_2, r_3}$. Отметим, что эти пересечения
		непустые, так как содержат $v$, и не двумерные, значит --- действительно прямые.
		
		Плоскость $\aff{l_1, l_2}$ является опорной для $C_1$ и содержит её вершину. Рассмотрим сечение
		$C_1$ плоскостью $\pi \parallel \aff{l_1, l_2}$. Обозначим $\pi \bigcap r_i = v_i$. Покажем, 
		что $v_1v_2v_3v_4$ --- параллелограмм. Действительно, $\pi \parallel l_1 \parallel \aff{r_1, 
		  r_2}$, поэтому $\pi \bigcap \aff{l_1, l_2} \parallel l_1$ и $v_1v_2 \parallel l_1$.
		Аналогично $v_3v_4 \parallel l_1$, $v_1v_4 \parallel l_2 \parallel v_2v_3$. Рассмотрим точки
		$v^*_1, v^*_2, v^*_3, v^*_4$, симметричные точкам $v_1, v_2, v_3, v_4$ в вершине $v$ веера.
		$v^*_1v^*_2v^*_3v^*_4$ --- также параллелограмм, его вершины лежат на одномерных рёбрах ячейки
		$C^*_1$ симметричной $C_1$. Четвёрка $v_1v_2v^*_4v^*_3$ --- также параллелограмм, так как
		$v_1v_2 = v_3v_4 = v^*_4v^*_3$ и $v_1v_2 \parallel v_3v_4 \parallel v^*_4v^*_3$. Эта четвёрка 
		также задаёт конус веера. Аналогично, остальные 3 трёхмерных конуса заданы параллелограммами 
		$v_2v_3v^*_1v^*_4$, $v_3v_4v^*_2v^*_1$ и $v_4v_1v^*_3v^*_2$.
	\hfill $\Box$ \par\bigskip

\subsection{Теорема о политопальности вееров, имеющих каноническую нормировку}
  Полный веер, по определению, задаёт нормальное локально-конечное разбиение пространства на
	полиэдры. Поэтому корректен вопрос о том, существует каноническая нормировка данного полного
	веера или нет.
	
	\begin{theorem}\label{theoremTightFansFansScalingPolytopality}~{
	  Каноническая нормировка полного заострённого веера существует тогда и только тогда, когда 
		этот веер политопален.
	}
	\end{theorem}
	
	\smallskip
	\noindent {\it Доказательство.} 
	  Докажем, что из политопальности полного заострённого веера следует, что для него существует 
		каноническая нормировка. 
		
		Пусть полный заострённый веер $\mathcal F$ представляется в виде веера 
		граней некоторого многогранника $P$. Пусть $v$ --- вершина $\mathcal F$. Тогда каждая $k$-грань 
		$f \prec \mathcal F$ имеет вид $f = \cone_v(F)$ для соответствующей $(k-1)$-грани $F \prec P$.
		
		Рассмотрим полярный многогранник $\polar{P}$ относительно центра $v \in \int{P}$. Каждая 
		гипергрань $h_i = h^{d-1}_i$ из веера $\mathcal F$ является пересечением ровно двух ячеек
		$c_{i, 1}, c_{i, 2} \in \mathcal F$. По определению, $h_i = \cone_v(H_i)$, $c_{i, 1} = 
		  \cone_v(C_{i, 1})$, $c_{i, 2} = \cone_v(C_{i, 2})$ для некоторого коребра	$H_i$ и гиперграней
		$C_{i,1}, C_{i, 2}$ в $P$. Гиперграням $C_{i,1}, C_{i, 2}$ соответствуют вершины $\polface{C_{
		  i, 1}}$ и $\polface{C_{i, 2}}$ в $\polar{P}$.
		
		Определим нормировку гиперграней веера $\mathcal F$. Нормировку $s(h_i)$ гиперграни $h_i \in 
		  \mathcal F$ положим равной $\left|\polface{C_{i,1}}, \polface{C_{i,2}}\right|$ --- длине 
		ребра $\polface{H_i}$. Покажем, что эта нормировка является канонической. 
		
		Пусть $h_1, \ldots, h_t$ --- все гиперграни из $\mathcal F$, которые содержат коребро	$f = 
  	  f^{d-2} \in \mathcal F$, занумерованные в одном из двух направлений циклического обхода 
		вокруг $f$ (см. определение \ref{definitionIntroCycling}). Обозначим ячейки $c_i \in \mathcal 
		  F$, которые содержат $f$, так, чтобы $h_i = c_i \bigcap c_{i+1}$. По определению считаем 
		$c_{t+1} := c_1$.
		
		Пусть $\n_{i, i+1}$ --- единичная нормаль к $h_i$, направленная от $c_i$ к $c_{i+1}$. По 
		лемме \ref{lemmaTightFansOrthogSubspace} имеем $\polface{C_i}\polface{C_{i+1}} = 
		  \polface{H_i} \bot \cone_v(H_i) = h_i$. Отсюда получаем, что $\n_{i, i+1} = 
			\dfrac{\overrightarrow{\polface{C_i}\polface{C_{i+1}}}}{\polface{C_i}\polface{C_{i+1}}}$.
		Тогда для кручения нормировки $s$ выполнено:
		$$\Delta_s(f) = \sum^t_{i=1}\n_{i, i+1}s(h_i) =  \sum^t_{i=1}
		  \dfrac{\overrightarrow{\polface{C_i}\polface{C_{i+1}}}}{\polface{C_i}\polface{C_{i+1}}}
			\cdot \polface{C_i}\polface{C_{i+1}} = \sum^t_{i=1}\overrightarrow{\polface{C_i}
			\polface{C_{i+1}}} = 0$$
	  Так как коребро $f$ выбрана произвольно, то, по определению, нормировка $s$ каноническая.
		
		Докажем, что из существования канонической нормировки для полного заострённого веера следует 
		его политопальность.
		
		Пусть $s$ --- каноническая нормировка веера $\mathcal F \subset \mathbb E^d$, и $c_1, \ldots,
		  c_m$ --- все ячейки этого веера. По теореме \cite[Теорема 2.3]{Gavriliuk_liftings} для каждой 
		ячейки $c_i$ существует женератриса $\mathcal G(\mathcal F, s, c_i)$ (поверхность в $\mathbb 
		  E^{d+1}$), которая является графиком (функции) женератрисы $G_i(x)$. По определению, $G_i(x) 
			= 0$ при $x \in c_i$ и $G_i(x) > 0$ при $x \notin c_i$. 
			
		Определим функцию $G(x) ~:=~ \sum\limits^m_{i=1} G_i(x)$. Тогда $G(x)$ также линейна на каждой 
		ячейке $c_i$. Кроме того, $G(v) = \sum G_i(v) = 0$ и $G(x) > 0$ для любой точки $x \neq v$ (так 
		как найдётся ячейка не содержащая $x$). Легко видеть, что $G(x)$ также женератриса для 
		$\mathcal F$. Обозначим женератрису-поверхность, соответствующую $G(x)$, через $\mathcal G$.
		
		Рассмотрим точечное множество $M = \epi{G} \bigcap \{x_{d+1} = 1\}$ --- пересечение надграфика 
		функции $G(x)$ с гиперплоскостью $\{x_{d+1} = 1\}$, параллельной $\mathbb E^d$. Надграфик 
		$\epi{G(x)}$ является $(d+1)$-мерным выпуклым конусом. Его вершина $v$ является единственной 
		точкой пересечения $\epi{G}$ с опорной гиперплоскостью $\{x_{d+1} = 0\}$. Отсюда следует, что 
		$M$ --- выпуклый ограниченный $d$-многогранник. 
		
		Каждой собственной грани $F \prec M$ соответствует собственная грань $\cone_v(F) \prec \epi{
		  G(x)}$. То есть $\cone_v(F)$ является гранью женератрисы $\mathcal G = \partial(\epi{G})$. 
		Проекция $\cone_v(F)$ на $\mathbb E^d$ --- это грань веера $\mathcal T$, которая является 
		конусом $\cone_v(\widehat{F})$, порождённым проекцией грани $F$ на $\mathbb E^d$. 
		Следовательно, $\mathcal F$ состоит из конусов, порождённых проекциями граней многогранника 
		$M$, и является веером граней для $\widehat{M}$ --- проекции многогранника $M$ на $\mathbb 
		  E^d$.
			
	\hfill $\Box$ \par\bigskip
	
	Следующий результат известен специалистам. В литературе он стандартно доказывается в специальных 
	случаях (см. например \cite{Fulton}). Из теорем \ref{theoremTightFansFansScalingPolytopality} 
	и \cite[Теорема 2.3]{Gavriliuk_liftings} этот результат следует в общем случае.
	
	\begin{corollary}\label{corollaryTightFansLiftingFanClassic}~{
	  Подъём полного заострённого веера существует тогда и только тогда, когда этот веер политопален.
	}
	\end{corollary}
	
	\begin{theorem}\label{theoremTightFansScalingExists}~{
	  У произвольного тесного веера размерности $d = 2$ или $d = 3$ существует хотя бы одна 
		каноническая нормировка.
	}
	\end{theorem}
	
	\smallskip
	\noindent {\it Доказательство.} 
	  Согласно предложению \ref{propositionTightFans2dPolytopal} все полные двумерные веера, в том
		числе и тесные, политопальны. По теореме \ref{theoremTightFansPolytopical}, все трёхмерные
		тесные веера также политопальны. Отсюда, согласно теореме \ref{theoremTightFansFansScalingPolytopality},
		все тесные веера размерности 2 и 3 имеют каноническую нормировку.
	\hfill $\Box$ \par\bigskip

  Теорему \ref{theoremTightFansFansScalingPolytopality} можно обобщить на случай произвольных
	разбиений пространства.
	
	\begin{theorem}\label{theoremTightFansPolytopality}~{
	  Пусть заданы нормальное локально-конечное разбиение $\mathcal T$ пространства $\mathbb E^d$ и 
		некоторая его $k$-грань $F \in \mathcal T$. Каноническая нормировка звезды $\mathop{St}(F)$ 
		существует тогда и только тогда, когда веер схождений в грани $F$ политопален.
	}
	\end{theorem}
	
	\smallskip
	\noindent {\it Доказательство.} 
	  Согласно определению \ref{definitionTightFansMeetingFan}, веер схождений в грани определён не
		однозначно, а с точностью до аффинного преобразования. Однако свойство политопальности
		сохраняется при невырожденных аффинных преобразованиях веера. Действительно, если веер 
		представляется в виде веера граней некоторого многогранника $M$, то образ этого веера при 
		невырожденном аффинном преобразовании $\varphi$ является веером граней многогранника 
		$\varphi(M)$. Следовательно, формулировка теоремы корректна.
		
		Рассмотрим реализацию веера схождений $\mathcal F_F = \mathcal F(F, \mathcal T, \mathcal S)$
		в $(d-k)$-мерной плоскости $\mathcal S$, которая ортогональна грани $F$ и проходит через
		внутреннюю точку $p \in \int{F}$. При сечении $(d-k)$-плоскостью $\mathcal S$ разбиение 
		$\mathcal T$ переходит в разбиение $\widetilde{\mathcal T}$ плоскости	$\mathcal S$, грань $F$ 
		переходит в точку $p = \widetilde{F}$, произвольная $n$-грань $G \succ F$ переходит в 
		$(n-k)$-грань $\widetilde{G}$, содержащую $p$ как вершину.
		
		Рассмотрим произвольную нормировку $s$ гиперграней из $\mathop{St}(F)$. Определим по $s$ 
		нормировку $s'$ гиперграней из веера схождений $\mathcal F_F$. Для произвольной гиперграни
		$\widetilde{H}$ положим $s'(\widetilde{H}) := s(H)$. Докажем равенство кручений нормировок $s$ 
		и $s'$.
		
		Рассмотрим произвольное коребро $\widetilde{R} = \widetilde{R}^{(d-k)-2}$ в $\mathcal F_F$. Его 
		содержат гиперграни $\widetilde{H}_1, \ldots, \widetilde{H}_t$. Единичные нормали $\mathbf{n}_i 
		  \in \mathbb E^d$ к гиперграням $H_i \in \mathcal T$ также ортогональны грани $F \prec H_i$. 
		По выбору, $F \bot \mathcal S$ и $\mathcal S$ является дополнительной плоскостью к $\aff{F}$. 
		Следовательно, $\mathbf{n}_i$ параллельны плоскости $\mathcal S$ и также являются единичными 
		нормалями к $\widetilde{H_i}$ в $\mathcal S$. Отсюда следует равенство кручений нормировок:
		$$\Delta_{s'}(\widetilde{R}) = \sum^t_{i=1}s'(\widetilde{H_i})\mathbf{n}_i = 
			  \sum^t_{i=1}s(H_i)\mathbf{n}_i = \Delta_s(R)$$
				
		Обратно, по нормировке $s'$, заданной для реализации веера схождений $\mathcal F_F$ в 
		$(d-k)$-плоскости $\mathcal S \bot F$, правило $s(H) := s'(\widetilde{H})$ задаёт нормировку 
		$s$ гиперграней звезды $St_{\mathcal T}(F)$. Очевидно, что в этом случае также выполнено 
		равенство нормировок $\Delta_s(R) = \Delta_{s'}(\widetilde{R})$ для соответствующих друг другу
		корёбер $R \in \mathcal T$ и $\widetilde{R} \in \mathcal F_F$.
		
		По определению, нормировка является канонической, когда все соответствующие кручения равны 0. 
		Таким образом, каноническая нормировка звезды $\mathop{St}(F)$ существует тогда и только тогда,
		когда существует каноническая нормировка полного заострённого веера $\mathcal F_F$. По теореме
		\ref{theoremTightFansFansScalingPolytopality}, каноническая нормировка для $\mathcal F_F$
		существует тогда и только тогда, когда веер $\mathcal F_F$ политопален.
	\hfill $\Box$ \par\bigskip

  \subsection{Веера и параллелоэдры}	
	Следующая теорема объясняет интерес к тесным веерам в нашей работе. Условия, наложенные на тесные
	веера, позволяют сыметировать геометрию схождений параллелоэдров без использования целочисленных 
	решёток. Это удаётся сделать благодаря соответствию понятий стандартной грани в случае 
	параллелоэдров и в случае тесных вееров.
	
	\begin{theorem}\label{theoremTightFansCorrespondence}~{
	  Пусть задан параллелоэдр $P$ и нормальное разбиение $\mathcal T_P$ пространства $\mathbb R^d$ на 
	  его параллельные копии. Веер схождений этого разбиения в произвольной грани является тесным.
	}
	\end{theorem}
	
	\smallskip
	\noindent {\it Доказательство.}
	  Рассмотрим ячейки разбиения, сходящиеся в грани $F$ разбиения $\mathcal T_P$ и пересечение любых 
	  двух $d$-ячеек $P_i, P_j$, содержащих $F$. По свойству разбиений на параллелоэдры, их пересечение 
	  будет некоторой стандартной гранью $H$ (в смысле теории параллоэдров), то есть центрально-симметричной 
	  гранью, центр которой является также центром симметрии всего разбиения $\mathcal T_P$. Но тогда, 
	  очевидно, разбиение также локально-симметрично в центре $H$, а значит, по лемме 
	  \ref{lemmaLocalSymm}, и в любой внутренней точке грани $H$.
	  
	  Если размерность $F$ равна нулю, то лемма уже доказана. Иначе размерность $F$ больше нуля и для 
	  построения ассоциированного веера необходимо рассмотреть сечение разбиения трансверсальным к $F$ 
	  подпространством $S$, проходящим через произвольную внутреннюю точку $F$. Такое сечение будет 
	  проходить и через внутреннюю точку грани $H$. Для этой точки есть окрестность, внутри которой 
	  разбиение $\mathcal T_P$ симметрично. Но тогда и сечение этой окрестности подпространством $S$ 
	  также является центрально-симметричным. То есть ассоциированный веер в $F$ также локально 
	  симметричен в грани $H$ (точнее, в её сечении подпространством $S$), а значит она стандартна. 
	  Симметричность сечений $C_i, C_j$, полученных из ячеек $P_i, P_j$, также очевидно следует из
	  симметричности до взятия сечения. Что и требовалось.
	\hfill $\Box$ \par\bigskip
	
	Из этой теоремы и результатов предыдущего параграфа вытекает важное для теории параллелоэдров 
	следствие.
	
	\begin{proposition}\label{propositionTightFansScalingExists}~{
	  Звёзды произвольных $(d-2)$- и $(d-3)$-граней разбиения на параллелоэдры имеют канонические
		нормировки (не менее одной).
	}
	\end{proposition}
	
	Это важное свойство используется для доказательства гипотезы Вороного для односвязной 
	$\delta$-поверхности \cite{GarGavMag}. Если аналогичное утверждение удастся доказать в
	размерности 4, то это может существенно упростить доказательство теоремы Ордина о
	гипотезе в случае 3-неразложимых параллелоэдров. 4-мерного случая достаточно, чтобы доказать, 
	что условия теоремы 8 в \cite{Ordin} выполнены, откуда следует справедливость гипотезы для
	3-неразложимого случая. Доказательства аналогичных утверждений в высших размерностях, по всей
	видимости, дадут продвижения в случае $k$-неразложимых параллелоэдров.
	
	Также следует отметить, что из теоремы \ref{theoremTightFansClassification} о 5 типах 
	трёхмерных тесных вееров и из теоремы \ref{theoremTightFansCorrespondence} следует классификация
	Делоне \cite{Delone} о 5 типах схождений параллелоэдров в гранях коразмерности 3. Для 
	доказательства этой классификации достаточно лишь отметить, что все найденные 5 типов тесных 
	вееров реализуются как веера схождений параллелоэдров. Это легко установить для разбиений уже на 
	трёхмерные параллелоэдры.


\begin{thebibliography}{0}
  \bibitem{Aurenhammer} F. Aurenhammer, ``A criterion for the Affine Equivalence of Cell Complexes 
	  in $\mathbb R$ and convex polyhedra in $\mathbb{R^{d+1}}$,  Discrete and Computational 
		Geometry, vol.2, p.49-64, 1987
		
  \bibitem{BuhshtaberPanov} V.M. Buchstaber, T.E. Panov, Toric Topology, AMS Math Surveys and 
	  Monographs, 204 (2015).
		
	\bibitem{Delone} Delaunay B.N., Sur l\'a partition reguli\'ere de l'espace a 4 dimension. Изв. 
	  АН СССР, (1929) No 1, 79-110, No 2, 147-164.
		
	\bibitem{Edelsbrunner} H. Edelsbrunner, R. Seidel, ``Voronoi Diagrams and Arrangements'', 
	  Discrete and Computational Geometry, v.1 n.1, p.25-44, 1986, doi: 10.1007/BF02187681
		
	\bibitem{Davis} C. Davis, The set of non-linearity of a convex piecewise-linear function, Scripta 
    Mathematica 24 (1959), 219-228
		
  \bibitem{Fulton} W. Fulton, Intorduction to Toric Varieties, Princeton University Press, 1993
	
	\bibitem{GarGavMag} A. Garber, A. Gavrilyuk, A. Magazinov, The Voronoi conjecture for 
    parallelohedra with simply connected $\delta$-surface. Discrete and Computational Geometry, 
		vol.53 iss.2, p.245-260, 2015, doi: 10.1007/s00454-014-9660-z
	
	\bibitem{Gavriliuk_liftings} 	А. А. Гаврилюк, “Геометрия подъемов разбиений евклидовых 
	  пространств”, Геометрия, топология и приложения, Сборник статей. К 70-летию со дня рождения 
		профессора Николая Петровича Долбилина, Тр. МИАН, 288, МАИК, М., 2015, 49–66
	
	\bibitem{Magazinov3Tight} А.Н. Магазинов, ``К теореме Делоне о классификации схождений 
	  параллелоэдров в гранях коразмерности 3'', Моделирование и анализ информ. систем, 20:4, (2015),
		71-80.
		
	\bibitem{McMullen} P. McMullen, Duality, sections and projections of certain Euclidean tiliings,
    Geom. Dedicata 49 (1994), 183-202
		
	\bibitem{Ordin} A. Ordine, Proof of the Voronoi conjecture on parallelotopes in a new special 
	  case. Queen's University, Kingston (2005)
		
	\bibitem{RyshRyb} Ryshkov S.S., Rybnikov Jr. K.A., The theory of quality translations with 
	  applications to tilings. Eur. J. Comb., 18(4):431-444, 1997.

  \bibitem{Ziegler} Циглер Г.М., Теория многогранников / Пер. с англ. под ред. Н.П. Долбилина. 
	  М.:МЦНМО, 2014.






 



  






\end{thebibliography}
\end{document}